\documentclass[a4paper,oneside,reqno]{amsart}

\usepackage{amsmath,amssymb,epsfig}

\newtheorem{thm}{Theorem}
\newtheorem{lem}{Lemma}

\newtheorem{coro}{Corollary}
\newtheorem{prop}{Proposition}

\newtheorem{quest}{Question}

\newcommand{\spanl}{\operatorname{span}}

\newcommand{\rev}{{\text{rev}}}

\newcommand{\pic}[1]{\raisebox{-0.8\baselineskip}%
{\epsfig{file=#1,height=2.0\baselineskip}}}

\title{The Lawrence--Krammer representation is unitary}
\author{Won Taek Song}
\keywords{Braid, unitary, representation;}

\email{cape@knot.kaist.ac.kr}

\address{Department of Mathematics, 
Korea Advanced Institute
of Science and Technology,
Taejon, 305-701, 
Korea}

\subjclass{Primary 20F36; Secondary 57M07.}

\begin{document}
\begin{abstract}
We show that the Lawrence--Krammer representation is unitary.
We explicitly present the non-singular matrix
representing the sesquilinear pairing  invariant under the action.
We show that reversing the orientation of a braid is 
equivalent to the transposition of its Lawrence--Krammer matrix followed
by a certain conjugation.
As corollaries it is shown that 
the characteristic polynomial of the Lawrence--Krammer matrix 
is invariant under substitution of its variables with their inverses
up to multiplication by units,
and is not a complete conjugacy invariant for braids.
\end{abstract}
\maketitle

\section{Introduction}
The Lawrence--Krammer representation 
$\mathcal{K}\colon B_n \to GL(\mathbf{Z}[t^{\pm1}, q^{\pm1}], n(n-1)/2)$
was first introduced
by Lawrence~\cite{MR92d:16020} 
and proved to be faithful by Bigelow~\cite{MR1815219} and 
Krammer~\cite{krammerbl}. 
As a braid invariant $\mathcal{K}$ is strong enough to distinguish
all braids. The characteristic polynomial $f_\beta 
= \det( z I_{n(n-1)/2} - \mathcal{K}(\beta))$, where $I_{n(n-1)/2}$ denotes
the $n(n-1)/2$ dimensional identity matrix, of the Lawrence--Krammer 
matrix~$\mathcal{K}(\beta)$ 
of a braid $\beta\in B_n$ appears to be rather good as a conjugacy
invariant of braids.
The author observed that $f$ does not detect the orientation reversal 
of strings of braids, hence is not a complete conjugacy invariant,
and that the polynomial $f_\beta$ has a symmetry
$f_\beta(z,t,q) \doteq f_\beta(z^{-1}, t^{-1}, q^{-1})$ just like
the Alexander polynomial of links (see Corollary~\ref{coro:incomplete} and 
\ref{coro:rev}).
Actually this paper arised in search of an explanation for 
these phenomena.

Let $D_n$ be an oriented disk with $n$ holes in the complex plane. 
The boundary $\partial D_n$ 
consists of 
$n$ puncture boundary components and 
the outer boundary component. 
Let $\overline{C}$ denote the space 
$\{ \{x, y\} \mid x,y\in D_n \}$
of all unordered pairs of points in $D_n$.
Let $ N(\operatorname{diag} \overline{C} ) $  be an
open regular neighborhood of  the subset 
$\operatorname{diag} \overline{C}  = \{ \{x,x\}\mid x\in D_n \}$.
Let $C = \overline{C} \setminus  N(\operatorname{diag} \overline{C} )$.
Then $C$ is a 4-manifold with boundary.


The fundamental group $\pi_1(C)$ is isomorphic to 
the subgroup $B_{n,2}$ of $B_{n+2}$ consisting of $(n+2)$-braids 
whose first $n$ strands go straight down  without winding 
and only the last two strands freely wind around.
Let $\phi\colon B_{n,2} \to \{t^a q^b \mid a,b\in \mathbf{Z}\}$
be the homomorphism to the abelian group of monomials
defined by $\phi(\gamma) = t^{s - 2l} q^l$ where $s$ is the exponent
sum of $\gamma$ in Artin generators $\sigma_i$
and $l$ is the linking number of the last two strands with the
first $n$ strands.
Note that $ \phi( \beta^{-1} \gamma \beta) = \phi(\gamma)$ 
for any $n$-braid $\beta \in B_n \subset B_{n,2}$.
Let $\tilde{C}$ be the covering space of $C$ whose fundamental
group is the kernel of $\phi$.

Let $\Lambda$ denote the ring $\mathbf{Z}[t^{\pm1}, q^{\pm1}]$
of two variable Laurent polynomials.
The homology group $H_2(\tilde{C})$ 
admits  a $\Lambda$--module structure.
It is a free $\Lambda$--module with rank $n(n-1)/2$ 
(see~\cite[Section~4]{MR1815219}).
The Lawrence--Krammer representation 
refers to the braid action on $H_2(\tilde{C})$
as $\Lambda$--module automorphisms. 

By the the Blanchfield duality theorem~\cite{MR19:53a} 
\cite[Appendix E]{MR97k:57011}, 
we have a non-degenerate sequilinear pairing
$\langle\ ,\ \rangle \colon BH_2(\tilde{C}, \partial \tilde{C}) \times 
BH_2(\tilde{C}) \to  \Lambda$
given by 
$$
\langle X, y \rangle = \sum_{a,b\in \mathbf{Z}} t^a q^b 
( X \cdot t^a q^b y)
$$
where $(\ \cdot \ )$ denotes the ordinary intersection number
and $BH$ denotes the torsion free part of a $\Lambda$--module $H$.
By identifying $BH_2(\tilde{C},\partial \tilde{C})$
with $\operatorname{Hom}(H_2(\tilde{C}), \Lambda)$,
we obtain a pairing 
$\langle\ ,\ \rangle \colon H_2(\tilde{C} ) \times 
H_2(\tilde{C}) \to  \Lambda$. Since the braid action on 
$H_2(\tilde{C})$
is induced from self-homeomorphisms of $\tilde{C}$,
clearly it should preserve the pairing $\langle\ ,\ \rangle$,
i.e., the Lawrence--Krammer representation is unitary. 
This could be a non-constructive and terse proof of
the following theorem.

\begin{thm}\label{thm:uni}
There exists a non-singular $n(n-1)/2 \times n(n-1)/2$ matrix $J$
over $\mathbf{Z}[t^{\pm1},q^{\pm1}]$ such that 
for each Lawrence--Krammer matrix $M = \mathcal{K}(\beta)$ of an
arbitrary $n$-braid $\beta\in B_n$
and its conjugate transpose $M^*$, 
the equality $ M J M^* = J $ holds.
\end{thm}
The symmetry property of the characteristic polynomial $f_\beta$
naturally follows.

In order to find the explicit matrix $J$ representing the pairing,
we need the concrete description of a base set of $H_2(\tilde{C})$ 
and its dual base in $BH_2(\tilde{C},\partial\tilde{C})$.
In~\cite{MR1804157,MR1815219} the \emph{forks} were used to express 
relative 2-cycles in $H_2(\tilde{C}, \partial \tilde{C})$.
Bigelow showed that given each relative 2-cycle 
$X \in H_2(\tilde{C}, \partial \tilde{C})$
defined by a fork, 
$(1-q)^2(1 + qt) X$ is an image of a 2-cycle 
from a closed surface immersed in 
$\tilde{C}$.
Using the set $\{v_{i,j}\mid 1\le i < j \le n  \}$ of
relative 2-cycles corresponding to \emph{standard forks} as a base,
he presented the explicit matrices for the Lawrence--Krammer representation.
Let $y_{i,j}$ be the 2-cycle in $H_2(\tilde{C})$ which maps to 
$(1-q)^2(1 + qt) v_{i,j}$.
We can define a $n(n-1)/2 \times n(n-1)/2$ 
matrix $J_1$ whose entries are given by
$\langle v_{k,l}, y_{i,j} \rangle$
for $1\le i < j \le n$ and $1 \le k < l \le n$.
Then for the Lawrence--Krammer matrix $M$ of an $n$-braid $\beta$
the equality $M J_1 M^*  = J_1$ follows from 
$ \langle v_{k,l} , y_{i,j} \rangle  =  
\langle \beta(v_{k,l}) , \beta(y_{i,j}) \rangle $.

The author did not succeed in obtaining a concrete description
of $y_{i,j}$ as a surface or in calculating the the matrix $J_1$.
The proof of Theorem~\ref{thm:uni} in this paper involves 
cumbersome matrix calculations on \emph{biforks}
defined in Section~\ref{sec:bifork}. 
The matrix $J$ of Theorem~\ref{thm:uni} is presented 
in Lemma~\ref{lem:multab} as a multiplication table of an algebra. 
Our proof has an advantage over the non-constructive one in that
the same method also works for
the proof of the following theorem, which relates the orientation
reversal of a braid to the matrix transposition.
As a corollory of this theorem we have $f_{\beta} = f_{\beta^\rev}$,
where $\rev \colon B_n \to B_n$ is the anti-automorphism
given by $\sigma_i^\rev = \sigma_i$.
\begin{thm}\label{thm:rev}
There exists a non-singular $n(n-1)/2 \times n(n-1)/2$ matrix $V$
over $\mathbf{Z}[t^{\pm1},q^{\pm1}]$ such that 
for each $n$-braid $\beta\in B_n$
the equality   
$  \mathcal{K}({\beta^{\rev}}) V  = V \mathcal{K}(\beta)^T  $
holds.
\end{thm}

In~\cite{MR1857374} it was shown that 
the Lawrence--Krammer representation is isomorphic to a summand of 
the representation on the Birman--Murakami--Wenzl (BMW) 
algebra~\cite{MR90g:57004,MR92k:17032,MR89c:57007}.
This implies that the Artin generators $\sigma_i$
satisfy certain cubic relations in the representation.
\begin{prop}\label{prop:cubic}
The Lawrence--Krammer representation satisfies the cubic relations
$$(\mathcal{K}(\sigma_i) - I_{n(n-1)/2}) 
(\mathcal{K}(\sigma_i)  + q I_{n(n-1)/2} ) 
(\mathcal{K}(\sigma_i) + t q^2 I_{n(n-1)/2} ) = 0 $$
for $i=1,\ldots,n-1$.
\end{prop}
The eigenvalue $-t q^2$ is associated to a 1--dimensional
eigenspace.
Therefore
$\sigma_i^{-1} 
(\sigma_i - 1) ( \sigma_i + q)  = \sigma_i - q \sigma_i^{-1} - (1-q) $
maps to a projector on to  the 1--dimensional subspace.
This fact  led the author to define 
the following skein relation of \emph{geometric biforks}.
$$
\pic{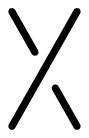} - q 
\pic{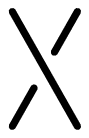} = (1-q)  \pic{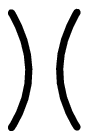} 
-  \pic{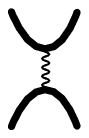}   
$$
Although
all the proofs in this paper were solely done by matrix calculations 
independent of the geometric biforks,
the definition of biforks and the main idea of the proofs came from
the geometric biforks.

In~\cite{MR99b:57015} it was shown that as a braid invariant
the Burau representation 
$\mathcal{B}\colon B_n \to GL(n-1,\mathbf{Z}[t^{\pm1}] ) $
is dominated by the 
finite type invariants. This is due to the fact 
that $\mathcal{B} ( \sigma_i  - \sigma_i^{-1} )(1) = 0$, 
i.e, at $t=1$
the Burau matrix is a finite type invariant of order 0.
Likewise the Lawrence--Krammer matrix
becomes a permutation matrix at $t=-1$ and $q=1$.
In Section~\ref{sec:fini} we show that
the matrix
$\frac{\partial^{k+l}}{ \partial t^k \partial q^l} \mathcal{K}(\beta)(-1,1)$
is of finite type of order $k+l$.
Therefore the Lawrence--Krammer representation is dominated
by the finite type invariants.
With the faithfulness of the representation together,
we obtain the following well-known fact.
\begin{prop}[\cite{MR89h:17030,MR96b:57004}]
The finite type invariants separate braids.
\end{prop}

The author is interested in the question whether the 
finite type conjugacy invariants of braids 
are strong enough to distinguish all the conjugacy classes of braids.
The characteristic polynomial $f_\beta$ of the Lawrence--Krammer
matrix is rather strong  as a conjugacy invariant, and
is dominated by the finite type invariants.
After observing that $f_\beta$ does not detect the 
orientation reversal of a braid, 
we raise the following question.
\begin{quest}
Is there a finite type conjugacy invariant $v$ of braids 
such that $v(\beta) \neq v(\beta^\rev)$ for some braid $\beta$?
\end{quest}

\section{Proofs of main results}\label{sec:bifork}
Let $\Lambda  = \mathbf{Z}[t^{\pm1},q^{\pm1}]$ be the ring of 
two variable Laurent polynomials in $t$ and $q$.
The Lawrence--Krammer module $V_n$ is a $n(n-1)/2$ dimensional free $\Lambda$--module 
generated by the set of standard forks $\{X_{i,j} \mid 1\le i < j \le n \}$.
The right action of $B_n$ on $V_n$ defined by~(\ref{eq:kram}) 
gives the faithful representation 
$\mathcal{K} \colon B_n 
\to  GL(n(n-1)/2, \Lambda)$~\cite{krammerbl,MR1815219}.
We follow the convention of~\cite{MR1815219} on the sign of variable~$t$.
\begin{equation}\label{eq:kram} 
 X_{i,j}\sigma_m = 
\begin{cases}
 X_{i,j} 
 \quad\hfill\text{if $m < i-1$ or $i<m <j-1$ or $j<m$    }  \\
 X_{i-1,j} 
 \quad\hfill 
 \text{if $m = i-1$ } \! \qquad  \\
 X_{i,j-1} 
 \quad\hfill 
 \text{if $i< m = j-1$ }  \\
-t q^2    X_{i,j} 
 \quad\hfill 
 \text{if $m=i=j-1$ }  \\
q  X_{i+1,j} + (1-q)   X_{i,j} 
 + t q (1-q)   X_{i,i+1} 
 \qquad \hfill \text{if $m=i < j-1$ }  \\
q  X_{i,j+1} + (1-q)   X_{i,j} 
 -  q (1-q)   X_{j,j+1} 
 \quad\hfill 
 \text{if $ m = j$ }  \ \quad \qquad
\end{cases}
\end{equation} 

Let $b\colon \{ (i,j) \mid 1 \le i < j \le n\} \to 
\{ \lambda \mid  1\le \lambda \le n(n-1)/2 \}$ be the bijection of index sets
defined by $b(i,j) = (2n-i)(i-1)/2 + (j-i)$. 
This is the lexicographic order on the double indices.
We identify $X_{i,j}$ with the 
$b(i,j)$-th standard base element
$Z_\lambda$, $\lambda = b(i,j)$,
which is the row vector whose $\lambda$-th entry  is $1$ and 
 uniquely non-zero.
$X_{i,j} \sigma_m$ denotes
the $b(i,j)$-th row of the matrix $\mathcal{K}(\sigma_m)$.

\begin{proof}[Proof of Proposition~\ref{prop:cubic}]
We need only to show the relation holds for $i=1$
since $\sigma_i$'s are all conjugate to each other. 

Firstly  it can be routinely verified for $j>2$  that 
$$q X_{1,j} \sigma_1^{-1} =  X_{1,j} - (1-q) X_{2,j}  
- ( 1- q^{-1}) X_{1,2}.$$
From this we calculate the images of $X_{1,j}$ and $X_{2,j}$
by $( \mathcal{K}(\sigma_1) - I_{n(n-1)/2} ) 
(\mathcal{K}(\sigma_i)  + q I_{n(n-1)/2} )$
as follows.
\begin{alignat*}{1}
X_{2,j} ( \mathcal{K}(&\sigma_1) - I_{n(n-1)/2} ) 
(\mathcal{K}(\sigma_i)  + q I_{n(n-1)/2} )  \\
&=   
X_{2,j} ( \mathcal{K}(\sigma_1) - q \mathcal{K}(\sigma_1^{-1})
 - (1-q) I_{n(n-1)/2} ) \mathcal{K}(\sigma_1)
\\
&= 
(X_{2,j} \sigma_1  - q X_{2,j} \sigma_1^{-1} - (1-q) X_{2,j} )
\mathcal{K}( \sigma_1) \\
&= 
( X_{1,j}  -   ( X_{1,j} - (1-q) X_{2,j}  
- ( 1- q^{-1}) X_{1,2}) - (1-q) X_{2,j}
) \mathcal{K}(\sigma_1) \\
&=  ( 1- q^{-1}) ( -t q^2)  X_{1,2}
\end{alignat*}
\begin{alignat*}{1}
X_{1,j} ( \mathcal{K}(&\sigma_1) - I_{n(n-1)/2} ) 
(\mathcal{K}(\sigma_i)  + q I_{n(n-1)/2} ) \hbox to5.2cm{ }  \\
&=   
X_{2,j} \sigma_1 ( \mathcal{K}(\sigma_1) - I_{n(n-1)/2} ) 
(\mathcal{K}(\sigma_i)  + q I_{n(n-1)/2} )  \\
&= 
X_{2,j}  ( \mathcal{K}(\sigma_1) - I_{n(n-1)/2} ) 
(\mathcal{K}(\sigma_i)  + q I_{n(n-1)/2} ) \mathcal{K}(\sigma_1)\\
&= 
( 1- q^{-1})  (-t q^2)  X_{1,2}\sigma_1 \\
&=  ( 1-q^{-1}) t^2 q^4   X_{1,2}
\end{alignat*}
We also have 
\begin{alignat*}{1}
X_{1,2}  ( \mathcal{K}(&\sigma_1) - I_{n(n-1)/2} ) 
(\mathcal{K}(\sigma_1)  + q I_{n(n-1)/2} )  = 
( -t q^2  - 1) ( -t q^2 + q) X_{1,2} \\
X_{i,j} ( \mathcal{K}(&\sigma_1) - I_{n(n-1)/2}) = 0
\text{\quad for $2 < i <j \le n$}
\end{alignat*}

The previous calculations exhibit that
$( \mathcal{K}(\sigma_1) - I_{n(n-1)/2} ) 
(\mathcal{K}(\sigma_1)  + q I_{n(n-1)/2} ) $
projects $V_n$ to the submodule generated by $X_{1,2}$.
The equality
$$(\mathcal{K}(\sigma_1) - I_{n(n-1)/2}) 
(\mathcal{K}(\sigma_1)  + q I_{n(n-1)/2} ) 
(\mathcal{K}(\sigma_1) + t q^2 I_{n(n-1)/2} ) = 0 $$
follows from 
$X_{1,2} (\mathcal{K}(\sigma_1) + t q^2 I_{n(n-1)/2} ) = 0 $.
\end{proof}

From the proof of Proposition~\ref{prop:cubic} we can see
that 
$  
 - (\mathcal{K}(\sigma_1) - I_{n(n-1)/2}) 
(\mathcal{K}(\sigma_1)  + q I_{n(n-1)/2}) 
\mathcal{K}(\sigma_1^{-1}) $
is an idempotent.
We denote this matrix by $X_{1,2}^* X_{1,2}^{\vphantom{*}}$.
$$X_{1,2}^* X^{\vphantom{*}}_{1,2} =  -  \mathcal{K}(\sigma_1) 
+ q \mathcal{K}(\sigma_1^{-1}) + (1-q) I_{n(n-1)/2}.$$
The first column of $X_{1,2}^* X^{\vphantom{*}}_{1,2} = Z_1^* Z^{\vphantom{*}}_1$ is 
the only non-zero column of $X_{1,2}^* X^{\vphantom{*}}_{1,2}$.
We denote by $X_{1,2}^*$ the  column vector
given by the first column of $X_{1,2}^* X^{\vphantom{*}}_{1,2}$
so that the expression $X_{1,2}^* X^{\vphantom{*}}_{1,2}$ makes sense also as 
a matrix multiplication.

Let $M(n(n-1)/2, \Lambda)$ be the set of matrices representing 
endomorphisms of $V_n$.
$M(n(n-1)/2, \Lambda)$ assumes a $B_n$-bimodule structure
by multiplications from the right side and from the left.
Let $\mathcal{M}_n$ be the $B_n$-bisubmodule of $M(n(n-1)/2, \Lambda)$
generated by the matrix $X_{1,2}^* X^{\vphantom{*}}_{1,2}$,
i.e., 
$$\mathcal{M}_n = 
\spanl_{\Lambda} \{\mathcal{K}(\beta) X_{1,2}^* X^{\vphantom{*}}_{1,2} \mathcal{K}(\gamma)
\in M(n(n-1)/2, \Lambda)  \mid \beta,\gamma \in B_n \}.$$

For $1\le i < j \le n$, 
let $A_{\pi(i,j)}$ denote the permutation braid
$(\sigma_{i-1}\cdots \sigma_1) 
(\sigma_{j-1} \cdots \sigma_2)$ \cite{MR96b:20052}. 
$A_{\pi(i,j)}$ induces the permutation $\pi(i,j)$ on $\{1,2,\ldots,n\}$ 
which maps $i$ to $1$, $j$ to $2$, and 
whose all  inversions involve either $i$ or $j$.
Let  $X_{i,j}^* X^{\vphantom{*}}_{k,l} = \mathcal{K}(A_{\pi(i,j)})
X_{1,2}^* X^{\vphantom{*}}_{1,2} \mathcal{K}( A_{\pi(k,l)}^{-1} ) $.
We call $X_{i,j}^* = \mathcal{K}( A_{\pi(i,j)} ) X_{1,2}^* $ 
a \emph{standard dual fork}.
We call $X_{i,j}^* X^{\vphantom{*}}_{k,l}$  a \emph{standard bifork}.
We denote $\mathcal{K}(\beta) X_{i,j}^* X^{\vphantom{*}}_{k,l} \mathcal{K}(\gamma)$
by $ \beta  X_{i,j}^* X^{\vphantom{*}}_{k,l} \gamma $
in emphasis of the right and left braid action on $\mathcal{M}_n$.

The following lemma is the multiplication table for 
standard biforks. 
\begin{lem}\label{lem:multab}
For $1\le i < j \le n$ and $1\le k < l \le n$,      
the $1\times 1$ matrix $X^{\vphantom{*}}_{k,l} X_{i,j}^*$ is given as follows.
$$
X^{\vphantom{*}}_{k,l} X_{i,j}^* = 
\begin{cases}
0 
  \kern12em
  \text{if $(i-k)(i-l)(j-k)(j-l) > 0$}  \label{lkj} \\
(-t^{-1} + q) ( q^{-1} + q t) 
\hfill{}
     \text{if $i=k < j=l$} \kern7.4em    \\
 t^{-1} (1 - q^{-1})
\hfill{}
    \text{if $i=k < l < j $ or  $k < i < j=l$} \kern0.2em \\
- (1 - q^{-1}) 
\hfill{}
    \text{if $i < j =  k < l$}\kern7.4em\\
  -t q (1-q) 
\hfill{}
    \text{if $i <  k < j = l$ or $i=k < j < l$}\kern0.2em\\
q (1-q)    
\hfill{}
   \quad\text{if $k < l= i < j$}\kern7.4em \\
- (1-q)^2 ( t^{-1} q^{-1} + 1) 
\hfill{}
    \text{if $k < i <  l < j$}\kern7.4em \\
(1-q)^2 ( q^{-1} + t) 
\hfill{}
   \text{if $i < k <  j < l$}\kern7.4em 
\end{cases}
$$
\end{lem}
\begin{proof}
This multiplication table can be verified
routinely from the definition of $X^{\vphantom{*}}_{i,j}$, $X_{k,l}^*$,
and the Lawrence--Krammer representation.
One may first calculate the table on $B_4$ (see~(\ref{eq:j4}))
in which every case of the double index correlations listed 
in the table occurs, and then use it for general $n$-braids.

We will not exhibit all the complicated calculations which could
be too distracting.
In Section~\ref{sec:geom} we introduce
geometric bifork algebra and calculate
the same multiplication table.
For each picture-based calculation line
in the proof of Theorem~\ref{thm:isom}
the same calculation on biforks can be done in parallel
by translating the isotopies to braid relations.
\end{proof}

\begin{lem}\label{lem:gen}
The set of standard biforks   $ \relax{B} = 
\{ X_{i,j}^* X^{\vphantom{*}}_{k,l} \mid 1 \le i < j \le n, 1 \le k < l \le n \}$
generates the algebra~$\mathcal{M}_n$ as a $\Lambda$--module.
\end{lem}
\begin{proof}
We need to show  that for all $\beta,\gamma \in B_n$
$\beta X_{1,2}^*  X^{\vphantom{*}}_{1,2} \gamma$ 
is a linear combination of standard biforks.
$X_{k,l} \gamma$ is a linear combination of standard forks
from the definition of the Lawrence--Krammer representation.
From the following formula we can see that 
$\beta X_{i,j}^*$ is also a linear combination of 
standard dual forks.

\begin{equation}\label{eq:invdual}
\sigma_m^{-1} X_{i,j}^*   = 
\begin{cases}
X_{i,j}^*  \kern6em  \text{if $m < i-1 $ or $i<m <j-1$ or $j<m$ }  \\
X_{i-1,j}^*  \hfill\text{if $m = i-1$ }\kern12em\kern-1pt  \\
X_{i,j-1}^*  \hfill\text{if $i< m = j-1$ } \kern10em\kern0.5pt \\
-t^{-1} q^{-2}   X_{i,j}^*  \quad \hfill\text{if $m=i=j-1$ }
   \kern10em\kern0.5pt   \\
q^{-1} X_{i+1,j}^*  + (1-q^{-1})  X_{i,j}^*  
 + t^{-1} q^{-1} (1-q^{-1})  X_{i,i+1}^* \\
\hfill   \text{if $m=i < j-1$} \kern10em\kern3.5pt  \\
q^{-1} X_{i,j+1}^*  + (1-q^{-1})  X_{i,j}^*  
 -  q^{-1} (1-q^{-1})  X_{j,j+1}^* \\
 \hfill\text{if $m = j$ } \kern13em\kern5pt  
\end{cases}
\end{equation}

\begin{enumerate}
\renewcommand{\theenumi}{\alph{enumi}}
\item
If  $m < i - 1$, then 
\begin{alignat*}{1}
\sigma_m^{-1} X_{i,j}^* X^{\vphantom{*}}_{1,2} &= 
     \sigma_m^{-1} A_{\pi(i,j)} X_{1,2}^* X^{\vphantom{*}}_{1,2} 
     = A_{\pi(i,j)} \sigma_{m+2}^{-1}  X_{1,2}^* X^{\vphantom{*}}_{1,2} \\ 
     &= A_{\pi(i,j)}  X_{1,2}^* X^{\vphantom{*}}_{1,2} \sigma_{m+2}^{-1} 
     = A_{\pi(i,j)}  X_{1,2}^* X^{\vphantom{*}}_{1,2} \\ 
     &= X_{i,j}^* X^{\vphantom{*}}_{1,2} 
\end{alignat*}
If  $ i < m < j-1$, then
$\sigma_m^{-1} A_{\pi(i,j)} = A_{\pi(i,j)} \sigma_{m+1}^{-1}$.  
If  $  j < m$, then 
$\sigma_m^{-1} A_{\pi(i,j)} = A_{\pi(i,j)} \sigma_{m}^{-1}$.  
For these two cases, the same calculation
works as for the case~$m < i - 1$.
\item
If $m = i-1$, then
$$\sigma_m^{-1} X_{i,j}^* X^{\vphantom{*}}_{1,2} = 
     \sigma_{i-1}^{-1} A_{\pi(i,j)} X_{1,2}^* X^{\vphantom{*}}_{1,2} 
     = A_{\pi(i-1,j)} X_{1,2}^* X^{\vphantom{*}}_{1,2} 
     = X_{i-1,j}^* X^{\vphantom{*}}_{1,2}. $$
\item
If $i < m = j-1$, then
 $$\sigma_m^{-1} X_{i,j}^* X^{\vphantom{*}}_{1,2} = 
     \sigma_{j-1}^{-1} A_{\pi(i,j)} X_{1,2}^* X^{\vphantom{*}}_{1,2} 
     = A_{\pi(i,j-1)} X_{1,2}^* X^{\vphantom{*}}_{1,2} 
     = X_{i,j-1}^* X^{\vphantom{*}}_{1,2}.$$
\item
If $ m= i = j-1$, then
\begin{alignat*}{1}
 \sigma_m^{-1} X_{i,j}^* X^{\vphantom{*}}_{1,2} &= 
     \sigma_{i}^{-1} A_{\pi(i,j)} X_{1,2}^* X^{\vphantom{*}}_{1,2} 
      = A_{\pi(i,j)} \sigma_1^{-1} X_{1,2}^* X^{\vphantom{*}}_{1,2} \\
      &= A_{\pi(i,j)} X_{1,2}^* X^{\vphantom{*}}_{1,2} \sigma_1^{-1} 
      = t^{-1} q^{-2} A_{\pi(i,j)} X_{1,2}^* X^{\vphantom{*}}_{1,2} \\
      &= t^{-1} q^{-2}  X_{i,j}^* X^{\vphantom{*}}_{1,2}.  
\end{alignat*}
\item
If $m=i < j-1$, then
\begin{alignat*}{1}
\sigma_m^{-1} X_{i,j}^* X^{\vphantom{*}}_{1,2} &= 
  \mathcal{K}(\sigma_i^{-1}) X_{i,j}^* X^{\vphantom{*}}_{1,2} \\ 
  &= \left( q^{-1} X_{i,i+1}^* X^{\vphantom{*}}_{i,i+1}   + q^{-1} \mathcal{K}(\sigma_i) 
+ (1-q^{-1}) I_{n(n-1)/2} \right) X_{i,j}^* X^{\vphantom{*}}_{1,2} \\
&=  q^{-1} X_{i,i+1}^* X^{\vphantom{*}}_{i,i+1} X_{i,j}^* X^{\vphantom{*}}_{1,2}   + 
q^{-1} \sigma_i A_{\pi(i,j)} X_{1,2}^* X^{\vphantom{*}}_{1,2} \\
&\ \hskip3cm  + (1-q^{-1}) X_{i,j}^* X^{\vphantom{*}}_{1,2}  \\
&=  q^{-1} X_{i,i+1}^* (   t^{-1} (1 - q^{-1})  ) X^{\vphantom{*}}_{1,2}   + 
q^{-1} A_{\pi(i+1,j)} X_{1,2}^* X^{\vphantom{*}}_{1,2} \\
&\hskip4cm + (1-q^{-1}) X_{i,j}^* X^{\vphantom{*}}_{1,2}  \\
&=   t^{-1} q^{-1}  (1 - q^{-1})   X_{i,i+1}^*  X^{\vphantom{*}}_{1,2}   + 
q^{-1}  X_{i+1,j}^* X^{\vphantom{*}}_{1,2} \\ 
& \hskip6cm  + (1-q^{-1}) X_{i,j}^* X^{\vphantom{*}}_{1,2}.  
\end{alignat*}
\item
If $m=j$, then
\begin{alignat*}{1}
\sigma_m^{-1} X_{i,j}^* X^{\vphantom{*}}_{1,2} &= 
  \mathcal{K}(\sigma_j^{-1}) X_{i,j}^* X^{\vphantom{*}}_{1,2} \\ 
  &= \left( q^{-1} X_{j,j+1}^* X^{\vphantom{*}}_{j,j+1}   + q^{-1} \mathcal{K}(\sigma_j) 
+ (1-q^{-1}) I_{n(n-1)/2} \right) X_{i,j}^* X^{\vphantom{*}}_{1,2} \\
&=  q^{-1} X_{j,j+1}^* X^{\vphantom{*}}_{j,j+1} X_{i,j}^* X^{\vphantom{*}}_{1,2}   + 
q^{-1} \sigma_j A_{\pi(i,j)} X_{1,2}^* X^{\vphantom{*}}_{1,2} \\
&\ \hskip3cm  + (1-q^{-1}) X_{i,j}^* X^{\vphantom{*}}_{1,2}  \\
&=  q^{-1} X_{j,j+1}^* (  - (1 - q^{-1})  ) X^{\vphantom{*}}_{1,2}   + 
q^{-1} A_{\pi(i,j+1)} X_{1,2}^* X^{\vphantom{*}}_{1,2} \\
&\hskip4cm + (1-q^{-1}) X_{i,j}^* X^{\vphantom{*}}_{1,2}  \\
&=  - q^{-1}  (1 - q^{-1})   X_{j,j+1}^*  X^{\vphantom{*}}_{1,2}   + 
q^{-1}  X_{i,j+1}^* X^{\vphantom{*}}_{1,2} \\ 
& \hskip6cm  + (1-q^{-1}) X_{i,j}^* X^{\vphantom{*}}_{1,2}.  
\end{alignat*}
\end{enumerate}
\end{proof}

The formula~(\ref{eq:invdual}) of $\sigma_m^{-1} X_{i,j}^*$ for the left action 
of $B_n$ on $\mathcal{M}_n$, presented in the proof of 
Lemma~\ref{lem:gen}, strikingly exhibits the duality
between forks and dual forks, and between the right action and 
the left action of braids.
If we replace in~(\ref{eq:invdual}) $\sigma_m^{-1}$ with $\sigma_m$, 
the left multiplication with the right one, 
$t$ with $t^{-1}$, 
$q$ with $q^{-1}$,
and $X_{i,j}^*$ with $X_{i,j}$, 
then we obtain exactly the same formula~(\ref{eq:kram}) 
used in the definition of the Lawrence--Krammer representation.
This observation interprets that if 
$Z_\lambda \sigma_m = \sum_\mu M_{\lambda\mu} Z_\mu$ then
\begin{equation}\label{eq:duality}
\sigma_m^{-1} Z_\lambda^* = \sum_\mu \overline{M_{\lambda\mu}} Z_\mu
\end{equation}
where $\overline{M(t,q)}$ denotes the matrix $M(t^{-1}, q^{-1})$.


\begin{lem}\label{lem:indep}
The set of standard biforks $\relax{B}$ is linearly independent.
\end{lem}
\begin{proof}
Let $J = \left( J_{\lambda\mu} \right)$ be the $n(n-1)/2 \times n(n-1)/2$ matrix
defined by 
$J_{\lambda\mu} =  Z_\lambda^{\vphantom{*}} Z_\mu^*$,
the multiplication table given in Lemma~\ref{lem:multab}.
At $q=1$, $J$ is a diagonal matrix with non-zero  diagonal entries, 
which implies $\det J \neq 0$ even for generic $q$.

Suppose that 
$\sum  a_{\lambda\mu} Z_\lambda^* Z^{\vphantom{*}}_\mu = 0$
for some $n(n-1)/2 \times n(n-1)/2$ matrix $A = \left( a_{\lambda\mu} \right)$.
Then for each $ 1\le \nu,o \le n(n-1)/2$
\begin{alignat*}{1}
0 =  Z^{\vphantom{*}}_\nu  
\left( \sum_{\lambda,\mu} a_{\lambda\mu} Z_\lambda^* 
Z^{\vphantom{*}}_\mu \right) Z_o^* 
 &= \sum_{\lambda,\mu} Z^{\vphantom{*}}_\nu Z_\lambda^*  
 a_{\lambda\mu} Z^{\vphantom{*}}_\mu Z_o^* \\
 &= \sum_{\lambda,\mu} J_{\nu\lambda} a_{\lambda\mu} J_{\mu o}, 
\end{alignat*}
 which implies $J A J = 0$.
Since $J$ is non-singular, $a_{\lambda\mu} =0$ for each 
$1\le \lambda,\mu \le n(n-1)/2$.
\end{proof}

Lemma~\ref{lem:indep} shows  in particular that 
the set of standard dual forks 
$\{Z_\lambda^* \mid 1 \le \lambda \le n(n-1)/2 \}$
is linearly independent.
Given a matrix $M(t,q) \in \mathcal{M}_n$, $M^*$ denotes 
the conjugate transpose  $M(t^{-1}, q^{-1})^T$. 
Here we assumed $t$ and~$q$ are evaluated at generic unit complex
numbers so that $t^{-1}$ and~$q^{-1}$ are complex conjugates of $t$ and~$q$.

\begin{proof}[Proof of Theorem~\ref{thm:uni}]
Let $J = \sum_{1\le i < j \le n}  X_{i,j}^* X^{\vphantom{*}}_{i,j}  
       =  \sum_{1 \le \lambda \le n(n-1)/2} 
       Z_\lambda^* Z^{\vphantom{*}}_\lambda$.

\begin{alignat*}{1}
J_{\lambda\mu}  &= Z_\lambda J Z_\mu^T   
       =  Z^{\vphantom{*}}_\lambda \left( \sum_{\nu }  
       Z_\nu^* Z^{\vphantom{*}}_\nu \right) Z_\mu^T \\
       &=  \sum_\nu Z^{\vphantom{*}}_\lambda Z_\nu^* 
       Z^{\vphantom{*}}_\nu Z_\mu^T 
        =  \sum_\nu Z^{\vphantom{*}}_\lambda Z_\nu^* \delta_{\nu\mu}\\
	&=  Z^{\vphantom{*}}_\lambda Z_\mu^* \label{eqn:jij}
\end{alignat*}
where $\delta_{\nu\nu} =1$ and $\delta_{\nu\mu} =0$ if $\nu\neq \mu$.
As shown above, $J$ is the same matrix given as the multiplication
table in the proof of Lemma~\ref{lem:indep}, where we 
showed that $J$ is non-singular.

It suffices to show the identity $ J M^* = M^{-1} J $
for the Lawrence--Krammer matrices $M = \mathcal{K}(\sigma_m)$ of 
the Artin generators $\sigma_m$ of $B_n$.

\begin{alignat*}{1}
M^{-1} J  &= \sum_\lambda  \sigma_m^{-1} Z_\lambda^* Z^{\vphantom{*}}_\lambda  
          = \sum_\lambda  \left(\sum_\mu Z_\mu^* \overline{M_{\lambda\mu}}\right) Z^{\vphantom{*}}_\lambda 
	    \quad \text{by (\ref{eq:duality})} \\
	  &= \sum_\mu  Z_\mu^* \left(\sum_\lambda  
	  Z^{\vphantom{*}}_\lambda \overline{M_{\lambda\mu}} \right)
	  = \sum_\mu  Z_\mu^*  
	  \left( Z^{\vphantom{*}}_\mu \overline{M}^T \right)  
	  = \left(\sum_\mu Z_\mu^* Z^{\vphantom{*}}_\mu \right) M^*\\
	  &= J M^*
\end{alignat*}
\end{proof}

\begin{coro}\label{coro:incomplete}
Let $f_\beta (z,t,q) = \det ( z I_{n(n-1)/2}  - \mathcal{K}(\beta)(t,q) )$
be  the characteristic polynomial of~$\mathcal{K}(\beta)$. 
Then  $f_\beta(z^{-1},  t^{-1} ,q^{-1})$ equals 
$f_\beta (z,t,q)$ up to  multiplication by units 
in~$\mathbf{Z}[t^{\pm1},q^{\pm1}]$.  
\end{coro}
\begin{proof}
By Theorem~\ref{thm:uni},  $ \mathcal{K}(\beta)^* 
= J^{-1} \mathcal{K}(\beta)^{-1} J$ as matrices over the quotient 
field~$\mathbf{Q}(t,q)$.
\begin{alignat*}{1}
f_\beta(z^{-1},  t^{-1},q^{-1}) 
    &= \det ( z^{-1} I_{n(n-1)/2}  - \mathcal{K}(\beta)^*  )\\
    &= \det ( z^{-1} I_{n(n-1)/2}  - J^{-1} \mathcal{K}(\beta)^{-1} J ) \\
    &=  \det ( z^{-1} I_{n(n-1)/2}  -  \mathcal{K}(\beta)^{-1}  ) \\
    &= (-z)^{-n(n-1)/2}  \det (\mathcal{K}(\beta)^{-1})  
        \det( z  I_{n(n-1)/2}  -  \mathcal{K}(\beta)  ) \\
    &= (-z)^{-n(n-1)/2}  \det (\mathcal{K}(\beta)^{-1})  f_\beta (z,t,q)
\end{alignat*}
\end{proof}

Let  $Y_{i,j} = X_{n+1-j, n+1-i} \Delta_n$ for $1\le i < j \le n$, 
where $\Delta_n$ denote the square root of
the full twist that generates the center of $B_n$.
For each $\sigma_m \in B_n$, we have 
$\Delta_n \sigma_m \Delta_n^{-1} = \sigma_{n-m}$.

\begin{lem}\label{lem:invkram}
With the base set $\{ Y_{i,j} \mid 1\le i < j \le n\}$, the 
Lawrence--Krammer representation is given by following formula.
\begin{equation}\label{eq:dualkram} 
 Y_{i,j}\sigma_m^{-1} = 
\begin{cases}
 Y_{i,j} \kern7em \text{if $m < i-1$ or $i<m <j-1$ or $j<m$    }  \\
 Y_{i-1,j} \hfill \text{if $m = i-1$ } \kern12em \kern-2pt  \\
 Y_{i,j-1} \hfill \text{if $i< m = j-1$ } \kern10em  \\
-t^{-1} q^{-2}    Y_{i,j} \hfill \text{if $m=i=j-1$ } \kern10em \\
q^{-1}  Y_{i+1,j} + (1-q^{-1})   Y_{i,j} 
 + t^{-1} q^{-1} (1-q^{-1})   Y_{i,i+1} \\
 \hfill \text{if $m=i < j-1$ } \kern10em \\
q^{-1}  Y_{i,j+1} + (1-q^{-1})   Y_{i,j} 
 -  q^{-1} (1-q^{-1})   Y_{j,j+1} \\
 \hfill \text{if $ m = j$ } \kern13em \kern4pt
\end{cases}
\end{equation} 
\end{lem}

\begin{proof}
From the definition~(\ref{eq:kram}) of the Lawrence--Krammer representation
the formula for $X_{i,j} \sigma_m^{-1}$ can be easily derived
as follows.
$$
X_{i,j}\sigma_m^{-1} = 
\left\{\hskip-0.5em \begin{minipage}[c]{.75\textwidth}
\vskip-.5\baselineskip
\begin{alignat}{1}
& X_{i,j}  \qquad\qquad \quad
  \text{if $m < i-1$ or $i<m <j-1$ or $j<m$ } \label{eq:inv1}  \\
& X_{i+1,j}
 \qquad\qquad  \text{if $m = i < j-1$ } \label{eq:inv2}  \\
& X_{i,j+1} \qquad\qquad   \text{if $   j  = m$ } \label{eq:inv3}  \\
& -t^{-1} q^{-2}    X_{i,j} \quad \!  \text{if $m=i=j-1$ } 
\label{eq:inv4} \\
& q^{-1}  X_{i-1,j} + (1-q^{-1})   X_{i,j} 
 -  q^{-1} (1-q^{-1})   X_{i-1,i} \label{eq:inv5}  \\
&\qquad \qquad \qquad  \quad \text{if $m = i-1  $ }\nonumber  \\
& q^{-1}  X_{i,j+1} + (1-q^{-1})   X_{i,j} 
 + t^{-1} q^{-1} (1-q^{-1})   X_{j,j+1}  \label{eq:inv6} \\
&\qquad \qquad \qquad \quad  \text{if $ i < m = j-1$} \nonumber
\end{alignat}
\end{minipage}
\right.
$$

The formula of this lemma is verified by the following 
routine calculations. 
\begin{enumerate}
\renewcommand{\theenumi}{\alph{enumi}}
\item 
If either $m < i-1$, $i<m <j-1$ or $j<m$,   
then $n+1-i <  n-m$, $ n+1 -j < n - m < (n+1-i) -1$, or
$n-m < (n+1-j)-1$ respectively.
Therefore 
\begin{alignat*}{1}
Y_{i,j}\sigma_m^{-1} &=  X_{n+1-j, n+1-i} \Delta_n \sigma_m^{-1} 
  = X_{n+1-j, n+1-i}  \sigma_{n-m}^{-1} \Delta_n \\
  &= X_{n+1-j, n+1-i}   \Delta_n 
  = Y_{i,j} \quad\text{by (\ref{eq:inv1}).}
\end{alignat*}
\item
If $m=i-1$, then $n-m = n+1-i$ so that by~(\ref{eq:inv3})
\begin{alignat*}{1}
Y_{i,j}\sigma_m^{-1}  &= X_{n+1-j, n+1-i}  \sigma_{n-m}^{-1} \Delta_n  \\
  &= X_{n+1-j, (n+1-i) + 1}\Delta_n  = Y_{i-1,j}.
\end{alignat*}
\item
If $i< m = j-1$, then $n-m =  n+1-j    <  (n+1-i) -1$
so that by~(\ref{eq:inv2})
\begin{alignat*}{1}
Y_{i,j}\sigma_m^{-1}  &= X_{n+1-j, n+1-i}  \sigma_{n-m}^{-1} \Delta_n  \\
  &= X_{(n+1-j) +1, n+1 -i} \Delta_n  = Y_{i,j-1}.
\end{alignat*}
\item
If $m=i=j-1$, then $n-m = n+1-j =  (n + 1 - i) -1$ 
so that by~(\ref{eq:inv4})
\begin{alignat*}{1}
Y_{i,j}\sigma_m^{-1}  &= X_{n+1-j, n+1-i}  \sigma_{n-m}^{-1} \Delta_n  \\
  &=  -t^{-1} q^{-2} X_{n+1-j, n+1-i} \Delta_n = -t^{-1} q^{-2} Y_{i,j}.
\end{alignat*}
\item
If $m=i < j-1 $, then $ n+1-j <   n-m = (n+1-i) -1$ 
so that by~(\ref{eq:inv6})
\begin{alignat*}{1}
Y_{i,j}\sigma_m^{-1}  &= X_{n+1-j, n+1-i}  \sigma_{n-m}^{-1} \Delta_n  \\
 & =  (q^{-1}  X_{n+1 - j ,(n+1-i)  +1} + (1-q^{-1})   X_{n+1-j, n+1-i} \\
 &\qquad \qquad \qquad   
  + t^{-1} q^{-1} (1-q^{-1})   X_{n+1-i,(n+1-i) + 1} )  \Delta_n  \\
 &= 
q^{-1}  Y_{i-1,j} + (1-q^{-1})   Y_{i,j} 
 + t^{-1} q^{-1} (1-q^{-1})   Y_{i-1,i}. 
\end{alignat*}
\item
If $m=j$, then $ n-m  =  (n+1-j)-1 $ 
so that by~(\ref{eq:inv5})
\begin{alignat*}{1}
Y_{i,j}\sigma_m^{-1}  &= X_{n+1-j, n+1-i}  \sigma_{n-m}^{-1} \Delta_n  \\
&= ( q^{-1}  X_{(n+1-j)-1,n+1-i} + (1-q^{-1})   X_{n+1-j,n+1-i}  \\
 &\qquad \qquad \qquad    
   - q^{-1} (1-q^{-1})   X_{(n+1-j)-1,n+1-j} ) \Delta_n \\
&=  q^{-1}  Y_{i,j+1} + (1-q^{-1})   Y_{i,j} 
 -  q^{-1} (1-q^{-1})   Y_{j,j+1}.
\end{alignat*}
\end{enumerate}
\end{proof}

For a word $W$, $W^\rev$ denotes the reverse word of $W$.
For a braid $\beta = W(\sigma_i) \in B_n$ written as a word 
in Artin  generators, we define $\beta^\rev = W^\rev(\sigma_i)$. 
In other words $\beta \mapsto \beta^\rev$ is the 
anti-isomorphism given by $\sigma_i^\rev = \sigma_i$.
Geometrically this equals reversing the orientations of 
the strings of a braid.

Observe in the formula~(\ref{eq:dualkram}) of 
Lemma~\ref{lem:invkram} that if one replaces 
$\sigma_m^{-1}$, $Y_{i,j}$, $t^{-1}$ and $q^{-1}$ with
$\sigma_m$, $X_{i,j}$, $t$ and $q$, then  
one obtains exactly the same formula~(\ref{eq:kram}). 

\begin{thm}\label{thm:base}
There exists an invertible $n(n-1)/2 \times n(n-1)/2$ matrix $R$
over $\mathbf{Z}[t^{\pm1},q^{\pm1}]$ such that 
for each $n$-braid $\beta\in B_n$,
the equality 
$$ \mathcal{K}(\beta^{-1})(t,q)  
= R^{-1} \mathcal{K}(\beta^\rev)(t^{-1}, q^{-1}) R  $$ 
holds.
\end{thm}
\begin{proof}
Let $W_{b(i,j)} = Y_{i,j}$ for $1\le i < j \le n$
and define the matrix $R$ by $R_{\lambda\mu} = (W_\lambda)_\mu$, 
the $\mu$-th entry
of $W_\lambda$ for $1\le \lambda,\mu \le n(n-1)/2$.

It suffices to show the equality for $M  = \mathcal{K}( \sigma_m)$.
The previous observation on the similarity between Lemma~\ref{lem:invkram}
and the definition of the Lawrence--Krammer representation 
interprets that $ W_\lambda M^{-1} = \sum_{\nu} 
\overline{ M_{\lambda \nu} } W_\nu $.
Taking the $\mu$-th entry of each side,  we have
\begin{alignat*}{1}
( W_\lambda M^{-1} )_\mu &= \sum_{\nu} \overline{ M_{\lambda\nu} } 
(W_\nu)_\mu \\
(R M^{-1})_{\lambda\mu} &=  \sum_{\nu} 
\overline{ M_{\lambda\nu} } R_{\nu\mu} \\ 
R M^{-1} &= \overline{M} R.
\end{alignat*}

From the definition of $Y_{i,j}$, $R \mathcal{K}(\Delta_n^{-1}) = P $
for some permutation matrix $P$. Hence $R$ is invertible.
\end{proof}

\begin{proof}[Proof of Theorem~\ref{thm:rev}]
Let $V = \overline{RJ}$.
Then we have
\begin{alignat*}{1}
\mathcal{K}(\beta^\rev) V &=   \mathcal{K}(\beta^\rev) \overline{RJ} 
  = \overline{ \overline{\mathcal{K}(\beta^\rev)} R J} \\
  &= \overline{ R \mathcal{K}(\beta^{-1}) J}  
    \qquad\text{by~Theorem~\ref{thm:base}} \\
  &= \overline{R  J \overline{\mathcal{K}(\beta)^T } }
    \qquad\text{by~Theorem~\ref{thm:uni}}  \\
  &= \overline{RJ} \mathcal{K}(\beta)^T 
  = V \mathcal{K}(\beta)^T.
\end{alignat*}
\end{proof}

Since transposition and conjugation do not alter 
the characteristic polynomial of a matrix,
we obtain the following corollary from Theorem~\ref{thm:rev}.
\begin{coro}\label{coro:rev}
$\mathcal{K}(\beta)$ and $\mathcal{K}({\beta^\rev})$ have the
same characteristic polynomial.
\end{coro}

\section{Explicit Matrices}
In this section we exhibit how to compute the matrices
$J$ and $V$ of the main theorems for low braid index.

Let  $J_4 = \sum_{1\le i < j \le 4} X_{i,j}^*  X^{\vphantom{*}}_{i,j} \in \mathcal{M}_4$.
$B_4$ is generated by the two elements $\sigma_1$ and 
$\delta_4  = \sigma_3\sigma_2\sigma_1$.
Note that $\sigma_2 = \delta_4^{-1} \sigma_1 \delta_4$ and 
$\sigma_3 = \delta_4^{-2} \sigma_1 \delta_4^2$.
For the proof of Theorem~\ref{thm:uni} for 4-braids, 
it is enough to verify the equalities 
$\mathcal{K}(\sigma_1) J_4  \mathcal{K}(\sigma_1)^*  = 
\mathcal{K}(\delta_4) J_4  \mathcal{K}(\delta_4)^* = J$.
The three matrices $J_4$,   $\mathcal{K}(\sigma_1)$  and 
$\mathcal{K}(\delta_4)$ can be explicitly written as follows:



\begin{alignat*}{1}
\mathcal{K}(\sigma_1)  &= 
\begin{bmatrix}
-t q^2 & 0& 0& 0& 0& 0 \\   tq (1 - q) & 1 - q& 0& q& 0& 0 \\  
     tq(1 - q)& 0& 1 - q& 0& q& 0 \\  0& 1& 0& 0& 0& 0 \\  
     0& 0& 1& 0& 0& 0 \\  0& 0& 0& 0& 0& 1
\end{bmatrix}  \\
\mathcal{K}(\delta_4)  &= 
\begin{bmatrix}
 0& 0& 0& q^2& 0& 0 \\ 0& 0& 0& 0& q^2& 0 \\ -tq^2& 0& 0& 0& 0& 0 \\ 
     0& 0& 0& 0& 0& q^2 \\ 0& -tq^2& 0& 0& 0& 0 \\ 0& 0& -tq^2& 0& 0& 0
\end{bmatrix}
\end{alignat*}

\begin{equation}\label{eq:j4}
J_4 =
\begin{bmatrix} 
 X_{1,2}^* &  X_{1,3}^* &  X_{1,4}^* &
 X_{2,3}^* &  X_{2,4}^* &  X_{3,4}^*  
\end{bmatrix}
\end{equation}
where the $6\times1$ matrices $X_{i,j}^*$ are given by
\begin{center}
\begin{minipage}{0.4\textwidth}
$$
X_{1,2}^* = 
\begin{bmatrix}
(-t^{-1} + q) ( q^{-1} + q t) \\
  -t q (1-q) \\
  -t q (1-q) \\
- (1 - q^{-1}) \\
- (1 - q^{-1}) \\
0 
\end{bmatrix}
$$
\end{minipage}
\hfill
and 
\hfill
\begin{minipage}{0.4\textwidth}
\begin{alignat*}{1}
X_{1,3}^* &= \mathcal{K}(\sigma_2) X_{1,2}^* \\
X_{1,4}^* &= \mathcal{K}(\sigma_3\sigma_2) X_{1,2}^* \\
X_{2,3}^* &= \mathcal{K}(\sigma_1\sigma_2) X_{1,2}^* \\
X_{2,4}^* &= \mathcal{K}(\sigma_1 \sigma_3 \sigma_2) X_{1,2}^* \\
X_{3,4}^* &= \mathcal{K}(\sigma_2\sigma_1 \sigma_3 \sigma_2) X_{1,2}^*. 
\end{alignat*}
\end{minipage}
\end{center}

\vskip8mm
The Lawrence--Krammer matrices of Artin generators $\sigma_1,\sigma_2$ of $B_3$
can written as follows:
$$
\mathcal{K}(\sigma_1) = 
\begin{bmatrix}
-t q^2 & 0 & 0 \\
 t q (1-q) & 1-q & q \\
0 & 1 & 0
\end{bmatrix}
\hfill
\qquad
\qquad
\qquad
\hfill
\mathcal{K}(\sigma_2) = 
\begin{bmatrix}
 1-q & q  & - q (1-q) \\
1 & 0 & 0 \\
0 & 0 & -t q^2 
\end{bmatrix}
$$
Let $Y_{1,2}^* = X_{1,2}^*$ be the $3\times 1$ matrix taken from 
the first column of $X_{1,2}^* X^{\vphantom{*}}_{1,2} \in \mathcal{M}_3$.
Let $V_3 = \begin{bmatrix} Y_{1,2}^* & Y_{1,3}^* &Y_{2,3}^* \end{bmatrix}$
where $Y_{1,3}^* = \mathcal{K}(\sigma_2^{-1}) Y_{1,2}^*$
and $Y_{2,3}^* = \mathcal{K}(\sigma_1^{-1}\sigma_2^{-1}) Y_{1,2}^*$.
The three columns of $V_3$ are explicitly written as follows:

\begin{alignat*}{1}
Y_{1,2}^* &= 
\begin{bmatrix}
(-t^{-1} + q) ( q^{-1} + q t) \\
-tq (1 - q) \\
- (1 - q^{-1})  
\end{bmatrix}  \\
Y_{1,3}^* &= 
\begin{bmatrix}
-t q (1-q)  \\
-(1 - q + q^2) ( -t - q^{-1} + t^{-1} q^{-3}) - q \\
 t^{-1} q^{-2} (1- q^{-1})
\end{bmatrix} \\
Y_{2,3}^* &= 
\begin{bmatrix}
- (1 - q^{-1}) \\ 
 t^{-1} q^{-2} (1- q^{-1}) \\
(-t^{-1} q^{-2} + q^{-1}) ( q^{-1} +  t q) 
\end{bmatrix} 
\end{alignat*}
Let $\beta = \sigma_1^2 \sigma_2 \sigma_1^{-2} \sigma_2^{-1}$ in $B_3$.
Then $\beta^\rev =  \sigma_2^{-1} \sigma_1^{-2} \sigma_2 \sigma_1^2$
is related to $\beta$ by a flype move, which 
changes the conjugacy class while preserving the link type of
a closed braid~\cite{MR94i:57005}.
One may check the equality 
$\mathcal{K}(\beta^\rev) V_3 = V_3 \mathcal{K}(\beta)^T$
and that $\mathcal{K}(\beta^\rev)$ shares with $\mathcal{K}(\beta)$
the same characteristic polynomial.

\section{The Burau representation}
In this section we review the Squier's result 
that the Burau representation is unitary.
The reduced Burau representation 
$\mathcal{B} \colon B_n \to GL_{n-1}(\mathbf{Z}[t^{\pm1}])$
is defined by these two $(n-1)\times(n-1)$ matrices:
\begin{alignat*}{1}
\mathcal{B}(\sigma_1) &= I_{n-1} + t e_{1,2}  - (1+t)e_{1,1} \\
\mathcal{B}(\delta_n ) &= - t^{n-1} e_{n-1,1} 
          + \sum_{1\le i \le n-2} ( -t^{i} e_{i,1}  + e_{i,i+1})
\end{alignat*}
where $e_{i,j}$ denotes the elementary 
matrix whose only non-zero entry is the $(i,j)$~entry with
value $1$, 
and $\delta_{n} = \sigma_{n-1}\sigma_{n-2}\cdots \sigma_1 \in B_n$.

\begin{thm}[Squier]\label{thm:squier}
There exists a nonsingular $(n-1)\times(n-1)$ matrix~$J_0$
over $\mathbf{Z}[t^{\pm1}]$ such that for each
$\beta$ in $B_n$ it follows that 
$\mathcal{B}(\beta)^* J_0 \mathcal{B}(\beta) = J_0$.
\end{thm}
The reduced Burau representation can be interpreted 
as the action on the homology group $H_1(\widetilde{D_n} ; 
\mathbf{Z}[t^{\pm1}])$ of the infinite cyclic cover $\widetilde{D_n}$
induced by the braid homeomorphism.
It is natural to expect the reduced Burau representation is 
unitary because homeomorphisms should preserve intersection forms.
We clarify this point in the following proof.
The original proof  in~\cite{MR85b:20056} 
was done by giving the matrix $J_0$ and directly evaluating the 
equality.
\begin{proof}
Consider the pairing
$\langle \ ,\ \rangle \colon 
H_1(\widetilde{D_n}, \partial \widetilde{D_n}) \times 
H_1(\widetilde{D_n} ) \to \mathbf{Z}[t^{\pm1}]$
defined by 
$$ \langle x , y \rangle  =  \sum_k t^k ( t^k x \cdot y) $$
where $(\ \cdot\ )$ denotes the usual algebraic intersection number.
The pairing $\langle\ ,\ \rangle$ is sesquilinear.
The equalities 
$ \langle t x , y \rangle = t^{-1}  \langle  x , y \rangle  $
and $ \langle  x , t y \rangle = t  \langle  x , y \rangle  $
follow from the definition.
If an automorphism $h_*$ of 
$H_1(\widetilde{D_n}, \partial \widetilde{D_n})$
is induced by a homeomorphism $h\colon \widetilde{D_n} \to \widetilde{D_n}$,
then $\langle h_*(x) ,h_*(y) \rangle = \langle  x , y \rangle$

We embed $D_n$ in the complex plane $\mathbf{C}$
so that the $i$-th hole is placed around the point $i\in \mathbf{C}$.
Let $x_i \in \pi_1(D_n )$ denote the standard generator
represented by the closed curve winding only around the $i$-th hole
once.
Let $y_i \in  H_1(\widetilde{D_n} ) $ be the cycle
which is the lift of $x_i x_{i+1}^{-1}\in \pi_1(D_n)$.
Then $\{ y_i\mid 1\le i \le n-1 \}$ is a base for 
the free $\mathbf{Z}[t^{\pm1}]$--module $H_1(\widetilde{D_n} )$.
Let $w_i \in  H_1(\widetilde{D_n}, \partial \widetilde{D_n})$
denote the lift of
the relative cycle connecting the $i$-th puncture boundary
to the $(i+1)$-st one by a straight segment, for $1\le i \le n-1$.
Note that $y_i$ maps to $(t-1) w_i$ by the inclusion
$ H_1(\widetilde{D_n} ) \to
H_1(\widetilde{D_n}, \partial \widetilde{D_n})$.

Let $M = \mathcal{B}(\beta)$ so that 
$\beta (w_i ) = \sum_{k} M_{ki} w_k$ and
$\beta (y_j ) = \sum_{l} M_{lj} y_l$. 
Let $(J_0)_{ij} = \langle w_i, y_j \rangle$.
Then the equality $ J_0 =M^* J_0 M $ follows as below:
\begin{alignat*}{1}
\langle w_i, y_j \rangle  &=  \langle \beta(w_i), \beta(y_j) \rangle
= \sum_{k,l} \langle  M_{ki} w_k, M_{lj} y_l \rangle  \\
&= \sum_{k,l} \overline{M_{ki}} \langle   w_k,  y_l \rangle M_{lj} 
\end{alignat*}

That $\det J_0 \neq 0 $ can be shown easily by  evaluating $J_0$ at $t=0$.
\end{proof}
Given a relative cycle  
$v \in H_1(\widetilde{D_n}, \partial \widetilde{D_n})$,
the action of $\sigma_i \in B_n$ on $w$  is determined by
the intersection number $\langle w_i, v \rangle$ as follows:
$$
\sigma_i(v) =  v  + \langle w_i, v \rangle  y_i
$$
as pointed in~\cite{MR2001j:20055}. 
From this we obtain
\begin{alignat*}{1}
\sigma_i (w_j) &= w_j + \langle w_i, w_j \rangle y_i \\
        &= w_j + \langle w_i, w_j \rangle (t-1) w_i \\
        &= w_j + \langle w_i, (t-1)w_j \rangle  w_i \\
        &= w_j + \langle w_i, y_j \rangle  w_i. 
\end{alignat*}
This equation implies that the $i$-th row of 
the matrix $\mathcal{B}(\sigma_i) - I_{n-1}$
is the  unique non-zero row with $\langle w_i, y_j \rangle$ as 
the $j$-th entry. 
Therefore we calculate the intersection pairing $J_0$ by
$$
J_0 = \sum_{1\le i \le n-1} 
  \left( \mathcal{B}(\sigma_i) - I_{n-1} \right).
$$

The analogous formula for the Lawrence--Krammer representation
obtained in the previous section is:
\begin{alignat*}{1}
J &= \sum_{1\le i < j \le n} X_{i,j}^* X^{\vphantom{*}}_{i,j} \\
 &= \sum_{1\le i < j \le n} 
 ( I_{n(n-1)/2} - \mathcal{K}(b_{i,j}) ) 
(I_{n(n-1)/2}  + q \mathcal{K}(b_{i,j}^{-1})) 
\end{alignat*}
where $b_{i,j}$ denotes the band generator 
$A_{\pi(i,j)} \sigma_1 A_{\pi(i,j)}^{-1}$.


\section{Geometric biforks}\label{sec:geom}
In this section, we define geometric biforks
and a skein algebra $\mathcal{T}_n$ generated by the geometric biforks.

\begin{figure}
\epsfig{file=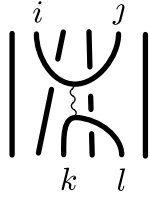}
\vskip-.5\baselineskip
\caption{a standard geometric bifork 
$\mathbf{x}_{i,j}^* \mathbf{x}_{k,l}$}\label{fig:bifork}
\end{figure}

Let $P = \{ p_i \in D^2  \mid 1\le i \le n \}$ be 
a set of $n$ distinct points in a disk $D^2$.
Then $P\times \{1, 0 \} \subset D^2 \times [0,1]$
is the set of $2n$ distinct points on the top and bottom of 
the solid cylinder $D^2\times [0,1]$.
A \emph{geometric $n$-braid} is a disjoint union of $n$ strings
in $D^2\times [0,1]$ having no local maxima or minima
with their end points fixed in $P\times \{1,0\}$.
A \emph{geometric bifork} is a disjoint union of $n$ strings
in $D^2\times [0,1]$ having exactly one local maximum 
and one local  minimum with one additionally attached string, which 
we call a \emph{handle}, connecting the maximum point to the minimum
point
without touching the other $n-2$ strings nor making 
a local extremum.
We distinguish the handle from the other strings 
by drawing it with a wavy line as in Figure~\ref{fig:bifork}.
The two strings to which the handle is attached are called 
\emph{tines}.
Two geometric biforks related by an isotopy 
which does not create a new local extremum are considered to be
the same.

We construct an arbitrary geometric bifork as follows.
Connect $p_1 \times \{1\}$ to $p_2 \times \{1\}$ 
by a string with exactly one local minimum in $D^2 \times [0,1]$ and
connect $p_1 \times \{0\}$ to $p_2 \times \{0\}$ 
by a string with exactly one local maximum. 
Then connect the local minimum to local maximum
by a straight wavy line.
For $3\le i \le n$, connect $p_i \times \{1\}$ to 
$p_i \times \{0\}$ by straight strings.
Now we obtained a simple geometric 
bifork $\mathbf{x}_{1,2}^* \mathbf{x}^{\vphantom{*}}_{1,2}$ which 
looks like a generator of the Birman--Murakami--Wenzl algebra 
except that it has a handle attached.
We may attach arbitrary braids $\beta$ and $\gamma$ 
to the top and the bottom of 
the geometric bifork to obtain a general geometric bifork 
$\beta \mathbf{x}_{1,2}^* \mathbf{x}^{\vphantom{*}}_{1,2} \gamma$ 
where $\beta, \gamma \in B_n$.
One can easily see that every geometric bifork can be written
as $\beta \mathbf{x}_{1,2}^* \mathbf{x}^{\vphantom{*}}_{1,2} \gamma$.

We define a $\Lambda$--algebra $\mathcal{T}_n$ generated by geometric biforks
with the following relations. 
\begin{alignat}{1}
\pic{sigma.eps} - q 
\pic{siginv.eps} &= (1-q)  \pic{triv.eps} 
-  \pic{idem.eps} \label{eq:f1}   \\
\pic{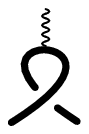} &= -t q^2  \pic{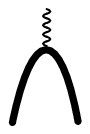} \label{eq:f2} \\ 
\pic{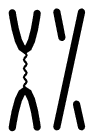} = \pic{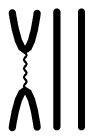}, & \qquad 
\pic{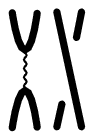} = \pic{idemtriv.eps}, \qquad
\pic{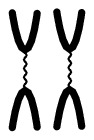} = 0 \label{eq:f3} \\
\pic{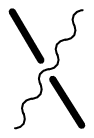} = q^2 \pic{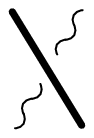}, &\qquad 
\pic{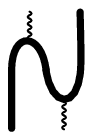} = q (1-q) \pic{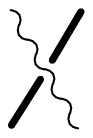},  
\qquad \pic{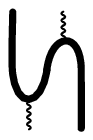} = q(1-q) \pic{negaphoton.eps} \label{eq:f4}
\end{alignat}
The multiplication in $\mathcal{T}_n$ is given by concatenation 
as in the braid groups.

The relations~(\ref{eq:f1},\ref{eq:f2}) come from the definition
that $X_{1,2}^* X^{\vphantom{*}}_{1,2} =  -  \mathcal{K}(\sigma_1) 
+ q \mathcal{K}(\sigma_1^{-1}) + (1-q) I_{n(n-1)/2}$
and $ X_{1,2} \sigma_1  = -t q^2 X_{1,2}$.
The relations~(\ref{eq:f3}) mean that 
$X_{1,2} \sigma_m = X_{1,2}$ for $m>2$.
The first relation in~(\ref{eq:f4}) reflects  that
$X_{1,2}^* X^{\vphantom{*}}_{1,2} \sigma_2\sigma_1^2 \sigma_2 
= q^2 X_{1,2}^* X^{\vphantom{*}}_{1,2}$
and the second one is from $X^{\vphantom{*}}_{1,2} X_{2,3}^* = q(1-q)$.

Let  $\mathbf{x}_{i,j}^* \mathbf{x}^{\vphantom{*}}_{k,l} = A_{\pi(i,j)} 
\mathbf{x}_{1,2}^* \mathbf{x}^{\vphantom{*}}_{1,2} A_{\pi(k,l)}^{-1} $.
We call $\mathbf{x}_{i,j}^* \mathbf{x}^{\vphantom{*}}_{k,l}$ a 
\emph{standard geometric bifork}. 
Figure~\ref{fig:bifork} shows a typical one.
Given an arbitrary geometric bifork $\mathbf{g}^*\mathbf{f}$, 
we can express $\mathbf{g}^*\mathbf{f}$ as a linear combination of 
other geometric biforks, which have less under-crossings of tines
than $\mathbf{g}^*\mathbf{f}$, by using
the relations~(\ref{eq:f1}--\ref{eq:f4}). 
Iterating this procedure, we express
$\mathbf{g}^*\mathbf{f}$ as a linear combination of standard
geometric biforks.
Therefore the set of standard geometric biforks generates $\mathcal{T}_n$
as a $\Lambda$--module.
One may check that 
the formulas~(\ref{eq:kram}) and (\ref{eq:invdual}) with 
$\mathbf{x}^{\vphantom{*}}_{i,j}$ and $\mathbf{x}_{k,l}^*$  in place of 
$X^{\vphantom{*}}_{i,j}$ and  $X_{k,l}^*$, also hold by applying 
the relations~(\ref{eq:f1}--\ref{eq:f4}).

The previous observation on the relations~(\ref{eq:f1}--\ref{eq:f4})
implies that there exists a surjective $\Lambda$--homomorphism
$\rho \colon  \mathcal{T}_n \to \mathcal{M}_n$, 
$\rho( \mathbf{x}_{i,j}^* \mathbf{x}_{k,l}^{\vphantom{*}} )
= X_{i,j}^* X_{k,l}^{\vphantom{*}}$.
\begin{thm}\label{thm:isom}
The geometric bifork algebra $\mathcal{T}_n$ 
is isomorphic to the matrix algebra $\mathcal{M}_n$ of biforks.
\end{thm}
\begin{proof}
By Lemma~\ref{lem:indep}, the algebra $\mathcal{M}_n$
is a free $\Lambda$--module with rank $(n(n-1)/2)^2$.
The fact that $\mathcal{T}_n$ is generated by
$( n(n-1)/2)^2$ many elements implies that
the surjective homomorphism 
$\rho \colon  \mathcal{T}_n \to \mathcal{M}_n$
is also injective.

In order to see $\rho$ is an algebra homomorphism, 
we need
$ (\mathbf{x}_{u,v}^*  \mathbf{x}_{k,l}^{\vphantom{*}} )
( \mathbf{x}_{i,j}^* \mathbf{x}^{\vphantom{*}}_{o,p} )
= 
(X_{k,l}^{\vphantom{*}} X_{i,j}^* ) \mathbf{x}_{u,v}^*  
 \mathbf{x}^{\vphantom{*}}_{o,p}$.
In the following we verify the equality
$  \mathbf{x}_{k,l}^{\vphantom{*}} \mathbf{x}_{i,j}^*  =
X_{k,l}^{\vphantom{*}} X_{i,j}^* $.
%
%
\begin{enumerate}
\renewcommand{\theenumi}{\alph{enumi}}
\item For $(i-k)(i-l)(j-k)(j-l) > 0$
$$
\pic{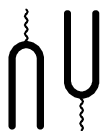} =
\pic{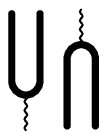} =
\pic{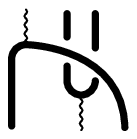} =
\pic{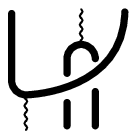} = 0
$$
\item For $i=k < j=l$
\begin{alignat*}{1}
\pic{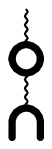} &=  - \pic{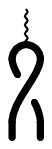}  + q \pic{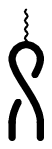} 
 + (1-q) \pic{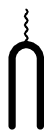} \\
 &= tq^2 \pic{nontw.eps} - t^{-1}q^{-1} \pic{nontw.eps} 
 + (1-q) \pic{nontw.eps}\\
 &= (-t^{-1} + q) (q^{-1} + q t) \pic{nontw.eps}
\end{alignat*}
\item For $i=k < l < j $ or  $k < i < j=l$
\begin{alignat*}{1}
\pic{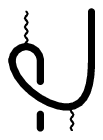}  &= - t^{-1} q^{-2} \pic{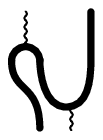} 
                   = -t^{-1} q^{-1} (1-q) \pic{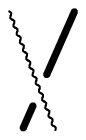} \\
\pic{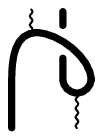}  &= - t^{-1} q^{-2} \pic{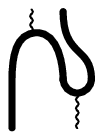} 
                   = -t^{-1} q^{-1} (1-q) \pic{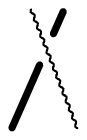} \\
\end{alignat*}
\item For $i < j =  k < l$
$$
\pic{x2x1.eps} = q (1-q) \pic{negaphoton.eps} = q^{-1} (1-q) \pic{photon.eps}
$$
\item For $i <  k < j = l$ or $i=k < j < l$
\begin{alignat*}{1}
\pic{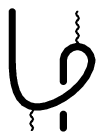}  &= - t q^2 \pic{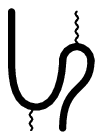}  
                   = -t q^3 (1-q) \pic{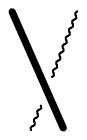}
		   = -t q (1-q) \pic{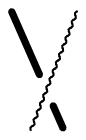} \\
\pic{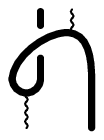}  &=   - t q^2 \pic{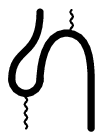} 
                   = -t q^3 (1-q) \pic{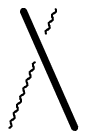} 
		   = -t q (1-q) \pic{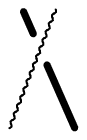} 
\end{alignat*}
\item For $k < l= i < j$
$$
\pic{x1x2.eps} = q (1-q) \pic{photon.eps} 
$$
\item For $k < i <  l < j$
\begin{alignat*}{1}
\pic{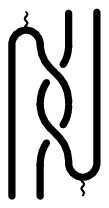} &=  q^{-1}\pic{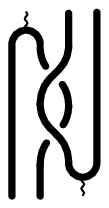} + (1-q^{-1}) \pic{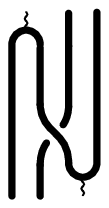}
  + q^{-1} \pic{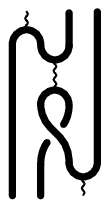} \\
  &= 0 + (1-q^{-1}) q (1-q) \pic{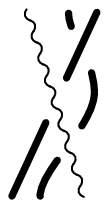}  
  - t^{-1}q^{-3} \pic{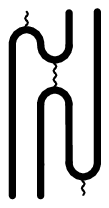} \\
  &= - (1-q)^2 \pic{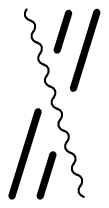} 
  - t^{-1} q^{-3} q^2 (1-q)^2 \pic{dirlink1.eps} \\
  &=  ( - (1-q)^2 - (1-q)^2 t^{-1} q^{-1}  ) \pic{dirlink1.eps} \\
  &=   - (1-q)^2(  t^{-1}  q^{-1} + 1 ) \pic{dirlink1.eps} 
\end{alignat*}
\item For $i < k <  j < l$
\begin{alignat*}{1}
\pic{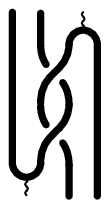} &=  q \pic{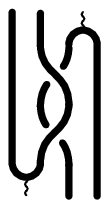} + (1-q) \pic{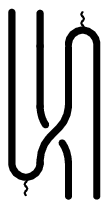}
  - \pic{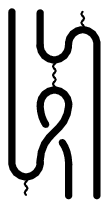} \\
  &= 0 + (1-q)  q (1-q) \pic{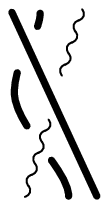}  
  +  t q^{2} \pic{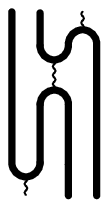} \\
  &= q^{-1}  (1-q)^2 \pic{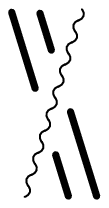} 
   + t q^{2} q^2 (1-q)^2 \pic{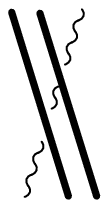} \\
  &=  ( q^{-1} (1-q)^2 +  t q^4 (1-q)^2 q^{-4}  ) \pic{dirlink2.eps}  \\
  &=    (1-q)^2(   q^{-1} + t ) \pic{dirlink2.eps} 
\end{alignat*}
\end{enumerate}

\end{proof}


\section{Dominance by finite type invariants}\label{sec:fini}
Let $I$ be the ideal of the integral group ring 
$\mathbf{Z}[B_n]$ generated by 
$\{ \sigma_i - \sigma_i^{-1} \mid 1\le i \le n-1 \}$.
Let $A$ be an abelian group.
If a $\mathbf{Z}$--module homomorphism $v \colon \mathbf{Z}[B_n] \to A$
vanishes on $I^{(k+1)}$,
we call $v$ a \emph{finite type invariant of order $k$}.
If $v(\gamma^{-1} \beta \gamma) = v(\beta)$ for each $\beta, \gamma \in B_n$,
we call $v$ a \emph{conjugacy invariant}.

\begin{thm}
For each $k,l\ge 0$, 
the $ n(n-1)/2 \times n(n-1)/2$ integral matrix invariant
$\mathbf{Z}[B_n] \to M( n(n-1)/2 , \mathbf{Z})$,
$$\beta \mapsto \frac{\partial^{k+l}}{\partial t^k \partial q^l} 
\mathcal{K}(\beta)(-1,1)$$
is a finite type invariant of order $k+l$.
\end{thm}
\begin{proof}
At $t=-1$ and $q=1$, 
it is easy to check that  
$\mathcal{K}( \sigma_m)(-1,1) 
=  \mathcal{K}(\sigma_m^{-1})(-1,1)$,
so that at $t=-1$ and $q=1$, $\mathcal{K}$ is a finite type invariant
of order $0$.
In the series expansion of
$\mathcal{K}( \sigma_m - \sigma_m^{-1}) 
= \sum_{k,l\ge 0} a_{kl} (t+1)^k (q-1)^l$ at $t=-1$ and $q=1$,
the lowest degree  is at least $1$.
In other words, 
$\mathcal{K}( W)$ has the lowest degree at least $1$ for each 
$W \in I^{(1)}$.

If $ W \in I^{(k+l+1)}$, then
the lowest degree in the series expansion 
of $\mathcal{K}( W)$ is at least $k+l+1$
since $W$ is a linear combination of  products of 
$(k+l+1)$ elements of $I^{(1)}$.
The theorem follows from 
$$\frac{\partial^{k+l}}{\partial t^k \partial q^l} 
(t+1)^i (q-1)^j |_{t=-1,q=1}=0 $$
for $i+j > k+l$.
\end{proof}


\providecommand{\bysame}{\leavevmode\hbox to3em{\hrulefill}\thinspace}

\end{document}